\documentclass[12pt]{article}

\pdfoutput=1

\usepackage{hyperref}
\usepackage{authblk} 
\usepackage{amsmath,amssymb}
\usepackage{graphicx}
\usepackage{color}

\usepackage{axisym-private} 
\date{}
\title{Numerical Strategies for Stroke Optimization of Axisymmetric
  Microswimmers\footnote{Preprint SISSA 33-2009-M.}}

\author[1]{Fran\c cois Alouges}
\author[2]{Antonio DeSimone} 
\author[2]{Luca Heltai\footnote{Corresponding Author. 
    Email: Luca Heltai \texttt{<luca.heltai@sissa.it>};  
    Tel.: +39 040-3787449; 
    Fax: +39 040-3787528}$^,$}

\affil[1]{CMAP UMR 7641, \'Ecole Polytechnique CNRS, \\
  Route de Saclay, 91128 Palaiseau Cedex - France}
\affil[2]{SISSA-International School for Advanced Studies, \\
  Via Bonomea 265, 34100 Trieste - Italy}

\begin{document}

\maketitle

\hrule
\section*{Abstract}


  We propose a computational
  method to solve optimal swimming problems, based on the boundary
  integral formulation of the hydrodynamic interaction between swimmer and surrounding fluid and direct
  constrained minimization of the energy consumed by the swimmer.

  We apply our method to axisymmetric model examples. We consider
  a  classical model  swimmer (the three-sphere swimmer of
  Golestanian et al.) as well as a
  novel axisymmetric swimmer inspired by the observation of
  biological micro-organisms.
\pagebreak

\section{Swimming at low Reynolds Numbers}
\label{sec:interpolation-v-g}

Swimming is the ability to advance within a fluid by performing cyclic
shape changes, in the absence of external propulsive forces. One of the
main difficulties of swimming at small scales is given by the
time-reversal property of Stokes flows, which describe hydrodynamics
at low Reynolds numbers.

The physical implication of this mathematical property is that simple
swimming strategies used in nature at larger scales, for example by
scallops, where the same shape change is executed forward and backward
at different velocities to achieve propulsion in one direction, do not
work for micro-, or nano-swimmers. At this scale, swimmers have to
undergo non-reciprocal deformations in order to achieve propulsion.
Also known as the Scallop-Theorem, this fact is discussed by Purcell
in~Ref.~\cite{Purcell-1977-a}, where a simple mechanical device
that can indeed swim at low Reynolds numbers, the three-link swimmer,
is proposed. A mathematical statement of the scallop theorem and its
proof can be found in Ref.~\cite{AlougesDeSimoneLefebvre-2009-a}.
More recently several other simple swimmers have been presented, for
example, in Golestanian and Najafi~\cite{NajafiGolestanian-2004-a} and
Avron et al.\cite{AvronGatKenneth-2004-a}

A mathematical approach to the problem of finding an optimal stroke
has been proposed by Alouges et
al. in Ref.~\cite{AlougesDeSimoneLefebvre-2007-a}, where it is shown
how to formulate and solve numerically the problem of finding optimal
strokes for low Reynolds number swimmers by focusing on the
three-sphere swimmer of Najafi and
Golestanian~\cite{NajafiGolestanian-2004-a} (a simple, yet
representative example).  The analysis carried out
in~Ref.~\cite{AlougesDeSimoneLefebvre-2007-a} shows how to address
quantitatively swimming as the problem of controlling shape in order
to produce a net displacement at the end of one stroke.  By casting it
in the language of control theory, the problem of swimming reduces to
the controllability of the system, and the search of optimal strokes
to an optimal control problem leading to the computation of suitable
sub-Riemannian geodesics.

The use of numerical algorithms to find optimal strokes can lead to
dramatic improvements. For the three-sphere swimmer, one can achieve
an increase of efficiency exceeding 300\% with respect to more naive
proposals.~\cite{AlougesDeSimoneLefebvre-2007-a} The simplicity of
three-sphere swimmer, which is a system with two degrees of freedom,
enables one to carry out the analysis in full detail. The study of
biologically relevant swimmers, however, requires more abstract
mathematical tools and more efficient numerical algorithms. The aim of
this paper is to introduce further numerical tools to overcome some of the
computational limits of the numerical method used
in~Ref.~\cite{AlougesDeSimoneLefebvre-2007-a}, opening up the
possibility to treat more and more complex swimmers such as the new
model swimmer described in section~\ref{sec:numerical-results}, which
is inspired by biological micro-swimmers like unicellular eukaryotic
algae.

For simplicity, we focus on the special case of axisymmetric swimmers,
which provides an interesting balance between complexity and
generality of the attainable results.  A complete theory to analyze
axisymmetric swimmers whose shape depends on finitely many parameters
and a general method to determine strokes of maximal efficiency is
presented in~Ref.~\cite{AlougesDeSimoneLefebvre-2008-a}. The main feature
qualifying this approach is the possibility of resolving hydrodynamic
interactions arising from the swimmer motion in their full complexity,
without being confined to asymptotic regimes (dilute limits for
assemblies of spheres, slender body hydrodynamics for one-dimensional
objects) in which explicit formulas for hydrodynamic forces are
available.

This work follows a similar line of thought. The main novelties
introduced here are the use of a boundary integral formulation for the
solution of the axisymmetric Stokes flow induced by the motion of the
swimmer, and the resolution of the optimal control problem via direct
constrained minimization of the functional giving the energy spent by
the swimmer. This should be contrasted with the resolution of the
hydrodynamic interactions via the finite element method that was used
in~Ref.~\cite{AlougesDeSimoneLefebvre-2007-a}
and~Ref.~\cite{AlougesDeSimoneLefebvre-2008-a}, and the resolution of the
Euler-Lagrange equation associated with the energy minimization
problem via an ODE scheme coupled with a shooting method to enforce
periodicity of the swimming stroke.

The accuracy that can be obtained with the new method is a few orders
of magnitude better than with the previous one, due to the possibility
of resolving the swimmer details at a much smaller scale at a
comparable computational cost.  In addition to the improved accuracy
granted by the use of boundary integrals, the direct solution of the
constrained minimization problem enlarges the class of optimization
problems we can consider.  In particular, it enables us 
to consider also inequality constraints on the shape variables. This
is often necessary in order to restrict the set of admissible shape
changes to those that are also biologically plausible.

The rest of this paper is organized as follows. In
section~\ref{sec:governing-equations} we derive the equations
characterizing optimal strokes, while in
section~\ref{sec:bound-integr-appr} we propose a strategy for their
numerical solution which is alternative to the ones presented
in~Ref.~\cite{AlougesDeSimoneLefebvre-2007-a}
and~Ref.~\cite{AlougesDeSimoneLefebvre-2008-a}. We validate our method by
comparing with previous results on the three-sphere swimmer, and by
calculating optimal strokes for a new model axisymmetric swimmer in
section~\ref{sec:numerical-results}.

\section{Governing Equations}
\label{sec:governing-equations}

We consider an axisymmetric micro-organism or micro-robot swimming at
low Reynolds numbers in a three dimensional fluid at rest at
infinity. No interaction with other obstacles is considered.

The swimmer occupies at time $t$ the (unknown) portion of space
denoted by $\Omega(t) \subset \Re^3$, with boundary $\partial
\Omega(t) = \Gamma(t)$, while the surrounding region
$\Re^3\backslash\Omega(t)$ is occupied by the fluid.

We are interested in finding the ``optimal stroke'' of period $T$,
directed along the axis of symmetry of the swimmer, and that will take
the swimmer from the position $\Omega_0$ to the position $\Omega(T) =
\Omega_0 + c\axis$, where $c/T$ is the average swimming velocity.

According to whether we want to optimize performance or energy
consumption, different optimality criteria can be adopted. In this
paper we select  the minimization of the
energy spent by a swimmer to swim at fixed (given) average velocity
$c/T$ as the optimality criterion.

\subsection{Exterior Hydrodynamic Interactions}
\label{sec:hydr-inter}

Independently of the optimality criterion that one chooses, any
swimmer needs to satisfy certain constraints which are required for
the problem to be well posed and for it to actually represent a
physically admissible swimmer.

An organism is said to be \emph{swimming} if it advances while its
shape changes periodically in time by following the hydrodynamics
of the surrounding fluid, in the absence of external forces or
torques.

More precisely, then, a  low Reynolds number swimmer satisfies the
following key properties:
\begin{itemize}
\item the velocity and force density at the surface of the swimmer
  are linked by the exterior Stokes equations satisfied by the surrounding fluid;
\item the shape of the swimmer changes periodically in time,
  driving the surrounding flow;
\item the force is generated internally by the swimmer, that is, the
  total external force and torque that act on the swimmer are zero.
\end{itemize}

On the exterior part of the swimmer, this is expressed by the
following system of equations, satisfied for each $t$ in $[0,T]$:
\begin{subequations}
 \label{eq:main-system}
  \begin{alignat}{2}
    & - \eta \Delta \u(\x,t) + \nabla p(\x,t) = -\nabla \cdot \vv
    \sigma = 0  && \mbox{in } \Re^3 \backslash \Omega(t) \label{eq:stokes-system} \\
    & \nabla\cdot \u(\x,t) =0 && \mbox{in } \Re^3 \backslash \Omega(t)\nonumber \\
    & \u(\x,t) = \v(\x,t) && \mbox{on } \Gamma(t)
    \nonumber \\
    & \lim_{|\x|\to\infty} |\u(\x, t)| = 0 &&
    \nonumber \\
    & \lim_{|\x|\to\infty} |p(\x, t)| = 0 && 
    \nonumber\\
    & \nonumber \\
    & \frac{\partial \X}{\partial t} = \v(\X,t) && \forall \X \in
    \Gamma(t)
    \label{eq:gamma-evolution}\\
    & \Gamma(T) = \vv R\Gamma(0) + \int_0^T \int_\Gamma\v \d \Gamma\d
    t
    = \vv R \Gamma(0) + c \axis \qquad&& \nonumber\\
    & \nonumber \\
    & \int_{\Gamma(t)} \vv \sigma \n = 0 && 
    \label{eq:internal-forces}\\
    & \int_{\Gamma(t)} \x \wedge \vv \sigma \n = 0 && 
    \nonumber
  \end{alignat}
\end{subequations}
where $\u$ and $p$ are the velocity and hydrodynamic pressure fields
in the exterior domain $\Re^3\backslash\Omega(t)$, $\eta$ is the
viscosity of the fluid, $\v$ is the velocity at the surface of the
swimmer, $\X$ is a material point on the surface $\Gamma(t)$ and
$\axis$ is the direction along which the swimmer is moving.

The set of Equations~\eqref{eq:stokes-system} describes the
conservation of momentum and mass in the Stokes fluid surrounding the
swimmer, with zero boundary conditions at infinity.

In Eq.~\eqref{eq:gamma-evolution} we grouped the \emph{no-slip}
boundary condition (each material particle of the surface of the
swimmer is following the evolution of the flow, i.e., no slip occurs
between the fluid and the swimmers) and the periodicity of the swimmer
shape, where $ c := \int_0^T\v\cdot\axis \d t$ is a shorthand notation
for the distance swan by the swimmer. Notice how in the general case
the configuration of the swimmer at time $T$ is a rigid motion of the
original configuration $\Gamma(0)$, where $\vv R$ indicates a rotation
and $c \axis$ is the translational vector.

Eq.~\eqref{eq:internal-forces} states that forces and torques are
generated internally by the swimmer. Here $\vv \sigma=\eta(\nabla\vv
\u+\nabla^T\vv \u)-p \vv I$ is the stress in the fluid, with $\vv I$
the identity.

In the axisymmetric case, the second of
Equations~\eqref{eq:internal-forces} is automatically satisfied, and
we will no longer refer to it. Moreover, due to the symmetry of the
system, no rotation is allowed and $\vv R$ is the identity matrix,
which we will omit from now on. Notice that this also implies that
both velocities and generated force densities are axisymmetric.

Knowledge of the velocity $\v(\x,t)$ on the boundary is enough to
ensure existence and uniqueness of a solution $\u(\x,t)$ in the entire
external domain, provided that $\Gamma(t)$ is maintained regular
enough throughout time: a boundary $\Gamma(t)$ of class $\mathcal{C}^1$
and boundary data $\v\in H^\frac12(\Gamma(t))$ guarantee existence and
uniqueness of the exterior solution $\u$ at time $t$, see, for
example~Ref.~\cite{DautrayLions-2000-a}.

\subsection{Global Hydrodynamic Interactions}
\label{sec:int-hydr-inter}

In Equations~\eqref{eq:stokes-system} we do not describe what
happens inside the domain $\Omega(t)$. While in reality micro
swimmers have a rather complex internal structure, we will make
two simplifying assumptions:
\begin{itemize}
\item all movement  capabilities of the organism are achieved through
  changes on the  surface $\Gamma(t)$ \emph{only};
\item the part of the swimmer enclosed by  $\Gamma(t)$ contains
fluid identical to the one in which the swimmer is immersed, hence
the same equations hold in the interior and in the exterior of
$\Omega(t)$.
\end{itemize}

Under these assumptions, a swimmer is fully defined if we either
assign a compatible movement of its boundary or by specifying the
distribution of stress jumps that the swimmer is sustaining on
$\Gamma(t)$.

The underlying principle we are exploiting here is that movement
is achieved via hydrodynamic interaction between the surface of
the swimmer and the surrounding fluid. Due to the complexity and
diversity of the internal structure of biological swimmers, we
don't attempt to model it explicitly. We do allow, however, for
the presence of a fluid inside $\Omega(t)$ because this is the
typical case for bio-swimmers. We assume for simplicity that the
viscosities of inner and outer fluids are the same because this
allows one to solve the global Stokes problem using only the
single layer potential. More precisely, the following relationship
holds (see, for example,~Ref.~\cite{Pozrikidis-1992-a}):
\begin{equation}
  \label{eq:integral-equation}
  8\pi\eta \u(\x) = - \int_\Gamma \vv{\tilde{\mathcal S}}(\x-\y) \jump{\vv
    \sigma(\y) \n(\y)} \d \Gamma(\y), \qquad \x \in \Re^3
\end{equation}
where $\vv{\tilde{\mathcal S}}$ is the three-dimensional
free-space Stokeslet  
\begin{equation}
  \label{eq:three-d-stokeslet}
  \vv{\tilde{\mathcal S}}(\vv r) := \frac{\mathbb{I}}{|\vv r|} +
  \frac{\vv r \otimes \vv r}{|\vv r|^3}
\end{equation}
and the fluid flow is solved for in the entire space $\Re^3$.

To simplify the notation a little, we introduce the definition of a
Neumann to Dirichlet $\ND$ and its inverse Dirichlet to Neumann map
$\DN$ as the operators that, for each stress jump distribution $\f =
\jump{\vv \sigma \vv n}$ on $\Gamma(t)$, return the velocity
distribution $\v$ generated on $\Gamma$ that satisfy
Eq.~\eqref{eq:integral-equation}, and \emph{vice versa}, i.e.:
\begin{multline}
  \label{eq:dirichlet-to-neumann-map}
  \ND\f  = \v  \Longleftrightarrow  \DN\v = \f  \Longleftrightarrow  \\ 
  8\pi\eta \v(\x) = - \int_\Gamma \vv{\tilde{ \mathcal S}}(\x-\y) \f(\y) \d
  \Gamma(\y), \, \x \in \Gamma .
\end{multline}
  
In the examples we consider, we always work with $\mathcal C^1$
boundaries $\Gamma$. Therefore Eq.~(\ref{eq:dirichlet-to-neumann-map})
makes sense for $\f \in H^{-\frac12}(\Gamma)$ and gives $\v\in
H^\frac12(\Gamma)$.

We remark here that both the Dirichlet to Neumann map and the Neumann
to Dirichlet map~\eqref{eq:dirichlet-to-neumann-map} are translation
invariant with respect to the configuration $\Gamma(t)$, i.e.,
\begin{equation}
  \label{eq:translation-invariance}
   \Gamma_1 = \Gamma_0+\vv d \Longrightarrow  N\!\!D_{\Gamma_0} \f =
    N\!\!D_{\Gamma_1} \f, \quad  D\!\!N_{\Gamma_0} \v =
    D\!\!N_{\Gamma_1} \v \qquad \forall \v, \quad \forall \f,
\end{equation}
for any constant translation vector $\vv d$, i.e., the same stress
jump applied to a translation of $\Gamma$ produces identical velocity
distributions on $\Gamma$ itself.

Given the above hypothesis, a biological swimmer can then be modeled
\emph{entirely} as a closed surface that is capable of transmitting
forces to the external (as well as to the \emph{internal}) fluid in
the form of stress jumps. The time evolution of the swimmer itself is
then given by transporting the boundary with the flow field $\u(\x,t)$
as in Eq.~\eqref{eq:gamma-evolution}, and the system of equations
satisfied by the unknown axisymmetric $\v$ and $\f$ becomes simply
\begin{subequations}
  \label{eq:main-system-reduced}
  \begin{alignat}{2}
    &\DNt \v = \f &&\label{eq:boundary-integral-equation-reduced}\\
    & \nonumber \\
   & \frac{\partial \X}{\partial t} = \v(\X, t) \qquad && \forall \X \in
   \Gamma(t)\label{eq:no-slip-bie}\\
  &\Gamma(T) = \Gamma(0) + \int_0^T \int_\Gamma \v \d\Gamma \d t =
  \Gamma(0) + c \axis
  \qquad&& \label{eq:periodicity-bie}\\
    & \nonumber \\
   & \int_{\Gamma(t)} \f= 0, &&
    \label{eq:boundary-integral-internal-forces}
  \end{alignat}
\end{subequations}
for $t \in [0,T]$, where
Eq.~\eqref{eq:boundary-integral-equation-reduced} effectively replaces
the set of Equations~\eqref{eq:stokes-system} and we already made the
simplifications derived from the symmetry of the system: no external
torque can be generated by axisymmetric force distributions, and the
only possible motion is a translation along the axis of symmetry.

More compactly, system~\eqref{eq:main-system-reduced} can be written
symbolically as
\begin{equation}
  \label{eq:compact-system}
  \begin{aligned}
    & \int_\Gamma \DN \left( \frac{\partial \Gamma }{\partial t}
    \right) \cdot \axis \d \Gamma = 0  \qquad && \forall t \in
    [0,T]\\
    &\Gamma(T) = \Gamma(0) + c \axis.
    \qquad
  \end{aligned}
\end{equation}

\begin{figure}[htb!]
  \centering
  \includegraphics[width=.8\textwidth]{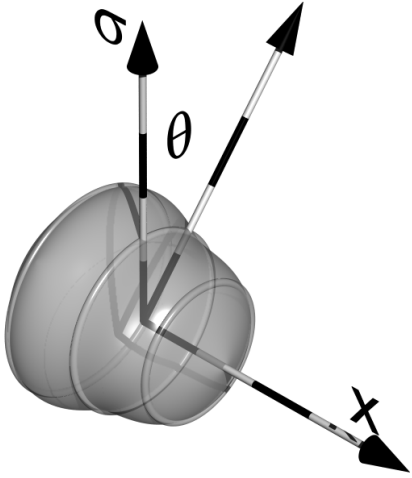}
  \caption{An example of the instantaneous configuration of an
    axisymmetric swimmer on the $(x,\sigma)-$plane. The three-dimensional
    configuration is obtained by a rotation around the $x$-axis of the
    function $\X(s)$.}
  \label{fig:configuration}
\end{figure}

Using cylindrical coordinates, the configuration of the swimmer, its
force density distribution and its velocity field on $\Gamma(t)$ can
all be expressed as a rotation around the axis of symmetry (i.e., the
$x$ axis) of one-parameter (vector valued) functions (see
Figure~\ref{fig:configuration} for an example configuration).

In particular we fix the domain $\Dom$ of the parameter $s$ to be the
interval $[0,1]$. The three-dimensional configuration $\Gamma$ is
the image of the function $\tilde\X$, given by
\begin{equation}
  \label{eq:axisymmetric-amoeba-configuration}
  \tilde\X(s,\theta) = \left(
    \begin{array}{c}
      x(s)\\
      \sigma(s) \cos(\theta) \\
      \sigma(s) \sin(\theta)
    \end{array}
  \right) = \Rot(\theta) \X(s), \quad \theta \in [0, 2\pi], \quad s \in [0,1],
\end{equation}
while the three-dimensional force density and the velocity fields are
given by
\begin{equation}
  \label{eq:axisymmetric-force}
  \tilde \f(s,\theta) = \left(
    \begin{array}{c}
      f_x(s) \\
      f_\sigma(s) \cos(\theta) \\
      f_\sigma(s) \sin(\theta)
    \end{array}
  \right) = \Rot(\theta) \f(s), \quad \theta \in [0, 2\pi], \quad s \in [0,1],
\end{equation}
and
\begin{equation}
  \label{eq:axisymmetric-velocity}
  \tilde \u(s,\theta) = \left(
    \begin{array}{c}
      u_x(s) \\
      u_\sigma(s) \cos(\theta) \\
      u_\sigma(s) \sin(\theta)
    \end{array}
  \right)  = \Rot(\theta) \u(s), \quad \theta \in [0, 2\pi], \quad s \in [0,1],
\end{equation}
where $\Rot(\theta)$ is defined as
\begin{equation}
  \Rot(\theta) = \left(
    \begin{array}{cc}
      1 & 0 \\
      0 & \cos(\theta) \\
      0 & \sin(\theta)
    \end{array}
  \right) \,.
\end{equation}

To avoid confusion, we will distinguish between three-dimensional
functions $\u(\x)$ and their axisymmetric counterparts $\u(s)$,
identified with their section on the $(x,\sigma)-$plane, by either
using symbols with a tilde on top (e.g., $\tilde {\mathcal S}$ for the
three-dimensional Stokeslet and $\mathcal S$ for the axisymmetric one)
or by explicitly specifying the domain variable ($\x$ for
three-dimensional functions or $s$ for their $(x,\sigma)$ section).

Our approach is to assign $\Gamma$ through a set of time dependent
scalars $\xi(t) \in V \subseteq \Re^N$, so that for each time $t$, the
curve $\X(\xi(t))(s)$ generates an admissible $\Gamma$. We focus on
cases where the boundary $\Gamma$ depends smoothly (e.g.,
analytically) on the parameters $\xi$, and remains of class
$\mathcal{C}^1$ for all time, $\X(s,t)$ is a non self-intersecting,
$\mathcal{C}^1$ curve on the $(x,\sigma)-$upper half plane ($\sigma
\geq 0$), and the surface obtained by its rotation around the $x$ axis
is the boundary $\Gamma(t)$ of a domain $\Omega(t)$ such that its
complement $\Re^3\setminus\Omega(t)$ is connected.

The admissible shapes we treat in this paper are collections of simple
closed curves in the $(x,\sigma)-$upper half plane, or open curves in
the same half plane which start and end vertically on the $x$ axis.

\subsection{Optimal Swimming as a Constrained Minimization Problem}
\label{sec:swimming-as-control-problem}

While system~\eqref{eq:main-system-reduced} describes the relationship
between the forces $\f$ that the swimmer applies to the surrounding
flow and its consequent evolution $\v$, it gives no information about
the \emph{optimality} of a given swimming stroke.

A natural optimality criterion consists on minimizing the energy
dissipated while swimming at a given average velocity $c/T$. In this
sense the optimal stroke is the one that minimizes
\begin{equation}
  \label{eq:minimization-problem}
  \EE(\v) := \int_0^T \int_{\Gamma(t)}\v\cdot \DN \v \d \Gamma \d t,
\end{equation}
where $\v$ satisfies the set of
Equations~(\ref{eq:main-system-reduced}), subject to the constraint
\begin{equation}
  \label{eq:constraint-displacement-velocity}
  \CC_c(\v) := \int_0^T \int_{\Gamma(t)}\v\cdot \axis \d \Gamma\d t -c = 0.
\end{equation}

In order to embed in the minimization problem also
Eq.~(\ref{eq:boundary-integral-internal-forces}), it is convenient to
express the entire problem in terms of shape changes rather than
absolute velocities.

We assume that the changes in the configuration $\Gamma(t)$ of the
swimmer happen \emph{only} through a set of $N+1$ scalar functions of
time which we identify as the shape variables $\xi_i(t), \, i=1\dots
N$, and the location $\varphi(t)$ on the axis of symmetry of a
distinguished point.

\begin{figure}[!htb]
  \centering
  \resizebox{\textwidth}{!}{\input{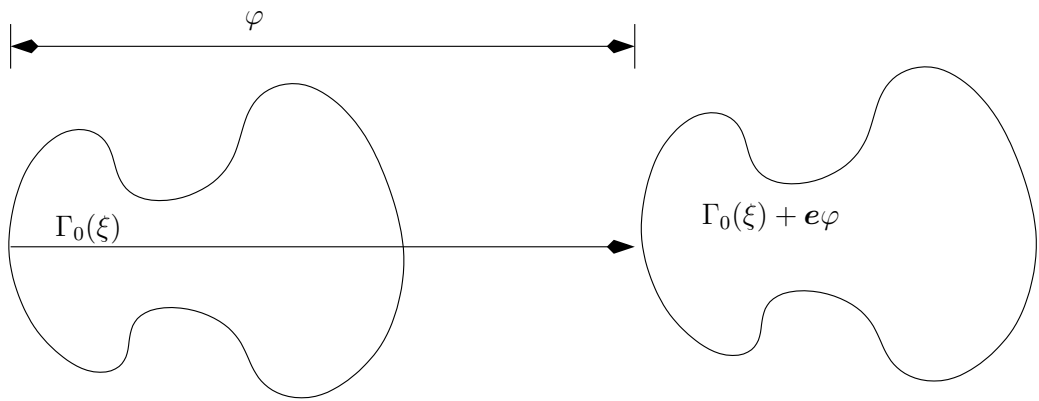}}
  \caption{Example of configuration $\Gamma(\xi, \varphi)$.}
  \label{fig:example-gamma}
\end{figure}

Moreover, the configuration $\Gamma(\xi, \varphi)$ is known for any
admissible shape $\xi$ and displacement $\varphi \in \Re$ and has
the form \begin{equation}
  \label{eq:gamma-of-xi-phi}
  \Gamma(\xi, \varphi) := \Gamma_0(\xi) + \axis \varphi, \qquad
  (\xi,\varphi) \in \Re^N\times\Re,
\end{equation}
as exemplified in Figure~\ref{fig:example-gamma}.

In particular, we assume that the configuration $\Gamma$ is the image
of a function $\X$ of the parameter $s \in \Dom$, which identifies the
location of each material point of $\Gamma$, i.e.,
\begin{equation}
  \label{eq:gamma-parametric}
  \X(\xi, \varphi, s) := \bar \X(\xi, s) + \axis \varphi, \quad \Gamma(\xi,
  \varphi) = \{ \X(\xi, \varphi, s),  \quad s\in \Dom \}\,.
\end{equation}

Using Equations~\eqref{eq:no-slip-bie} and~\eqref{eq:gamma-parametric},
the velocity of the swimmer can then be expressed in terms of rate of
shape changes $\dot \xi$ and in terms of translational velocity $\dot
\varphi$ as
\begin{equation}
  \label{eq:definition-velocity}
  \begin{split}
    \v(\X(\xi(t), \varphi(t), s) ) 
    =   &\frac{\partial \X(\xi, \varphi, s)}{\partial \xi_i} \dot \xi_i(t) +
    \frac{\partial \X(\xi, \varphi, s)}{\partial \varphi} \dot \varphi(t) \\
    =  &
    \frac{\partial \bar \X(\xi, s)}{\partial \xi_i} \dot
    \xi_i(t) + \axis \dot \varphi(t) \\
    = &
    \basexi_i(s) \dot\xi_i(t) + \axis \dot \varphi(t), \qquad \forall s\text{ in }
    \Dom, \quad \forall t \in [0,T],
\end{split}
\end{equation}
where summation is implied on repeated indices. Here $\axis$ is a unit
vector directed along the axis of symmetry, while $\basexi_i$ is the
change in the configuration $\Gamma(t)$ at the position $\X(\xi,
\varphi, s)$ due to a variation of the shape parameter $\xi_i$.

Using this formalism, the periodicity
condition~\eqref{eq:periodicity-bie} becomes simply
\begin{equation*}
 \xi(0) = \xi(T).
\end{equation*}

Eq.~(\ref{eq:boundary-integral-internal-forces}) implies that
\begin{equation}
  \label{eq:no-forces-condition}
  \int_\Gamma \DN \v \cdot \axis \d \Gamma = 0
\end{equation}
which can be rewritten using Eq.~\eqref{eq:definition-velocity}
as
\begin{equation}
  \label{eq:matrix-condition}
   \int_\Gamma \DN(\basexi_i \dot\xi_i + \axis \dot \varphi) \cdot \axis
   \quad  = \quad
   \NN(\xi) \cdot \dot \xi + \MM(\xi)\, \dot \varphi
   \quad  = \quad 0,
\end{equation}
where
\begin{equation}
  \label{eq:def-Phi-Xi}
  \begin{split}
    &  \NN _{i} (\xi) :=  \int_\Gamma \DN \basexi_i \cdot \axis
    \d \Gamma \\
    & \MM(\xi) := \int_\Gamma \DN \axis \cdot
    \axis \d \Gamma.
  \end{split}
\end{equation}

It is now possible to write $\dot \varphi$ as a function of $\dot \xi$
\begin{equation}
  \label{eq:definition-V}
  \dot \varphi = \VV(\xi) \cdot \dot \xi = - \frac{\NN(\xi)}{\MM(\xi)}\cdot \dot \xi,
\end{equation}
which allows us to express the energy dissipated by a given stroke
$\xi$ as
\begin{equation}
  \label{eq:energy-dissipated}
  \mathcal E(\xi) := \int_0^T \GG(\xi) \dot \xi \cdot  \dot \xi \d t,
\end{equation}
where
\begin{equation}
  \label{eq:def-G-w}
  \GG_{ij}(\xi) := \int_\Gamma \DN \basev_i \cdot \basev_j \d \Gamma, \qquad
  \basev_i := \basexi_i + \VV_i(\xi) \axis.
\end{equation}

A natural space where we can look for solutions $\xi$ of our problem
is then given by the Sobolev space
\begin{equation}
  \label{eq:solution-space}
  \Vs := \{ \xi \in H^1(0,T)^N, \text{ s.t. } \xi(0) = \xi(T) \}
\end{equation}
of periodic $H^1$ functions of $N$ components. 
Optimal swimming is then given by the solution $\bar \xi \in \Vs$ of
\begin{subequations}
  \label{eq:optimal-swimming}
  \begin{align}
    \min_{\xi \in \Vs}\qquad & \int_0^T \GG(\xi)\dot \xi \cdot \dot
    \xi \d t =: \EE(\xi) \label{eq:opt-minimality}\\
    \text{subject to }\qquad &\int_0^T\VV(\xi)\cdot\dot \xi \d t- c
    = : \CC_c(\xi) = 0 \label{eq:opt-feasibility}.
  \end{align}
\end{subequations}

In principle, additional constraints could be added to
Eq.~(\ref{eq:opt-feasibility}) to fix, for example, that the
swimmer initial and final shape are given, or that the shape never
reaches out of given physical bounds (due, for example, to
construction constraints in artificial swimmers or physiological
constraints in biological swimmers):
\begin{subequations}
  \begin{align}
    \label{eq:alternative-initial-condition}
    \xi(0) & = \xi(T) =\xi_0 \\
    \label{eq:shape-bounds}
    \xi_L & \leq \xi(t) \leq \xi_U & \forall t \in [0,T].
  \end{align}
\end{subequations}
In ~Ref.~\cite{AlougesDeSimoneLefebvre-2007-a}
and~Ref.~\cite{AlougesDeSimoneLefebvre-2008-a}, 
problem~(\ref{eq:optimal-swimming}) with the addition of
constraint~(\ref{eq:alternative-initial-condition}) 
has been considered for some special model swimmers.
In this case, the space of admissible shapes is 
\begin{equation}
  \label{eq:solution-space-withbc}
  \Vs_0 := \{ \xi \in H^1(0,T)^N, \text{ s.t. } \xi(0) = \xi(T)=\xi_0  \}
\end{equation}
and the following hypothesis guarantees the existence of a solution to
problem~(\ref{eq:optimal-swimming}):
\begin{equation}
\label{controllabilityhyp}
\mbox{There exists }\xi \in \Vs_0\mbox{ which satisfies }(\ref{eq:opt-feasibility})
\mbox{ and such that }\EE(\xi)\mbox{ is finite},
\end{equation}
as shown below.
\begin{proposition}
  Assume $\GG$ and $\VV$ are continuous functions and that there exist
  positive constants $\alpha, \beta$ such that
  \begin{equation}
    \exists \alpha>0, \beta>0,\text{ such that },\,\,
    \beta|\eta|^2\geq \eta\cdot \GG(\xi)\eta\geq \alpha |\eta|^2\quad\forall \eta \in \Re^N,
    \label{uniformcoercivity}
  \end{equation}
uniformly for admissible shapes $\xi\in \Vs_0$.
Then there exists a solution to problem (\ref{eq:optimal-swimming}).
\end{proposition}

\begin{proof}
  The proof is a straightforward application of the Direct Method of
  the Calculus of Variations. Indeed, under the given hypotheses, the
  functional to minimize is coercive and lower semicontinuous with
  respect to weak convergence in $(H^1(0,T))^N$. Moreover the
  constraint is continuous with respect to this topology. Condition
  (\ref{controllabilityhyp}) ensures that the set over which the
  minimization is done is not empty.
\end{proof}

\begin{remark}
  Condition (\ref{controllabilityhyp}) is nontrivial. It 
  requires that 
  the swimmer can indeed cover the prescribed distance, using finite energy.  As
  explained in Ref.~\cite{AlougesDeSimoneLefebvre-2007-a} and
  Ref.~\cite{AlougesDeSimoneLefebvre-2008-a}, this is a
  controllability condition equivalent to asking that the constraint
  (\ref{eq:opt-feasibility}) is non-holonomic.
\end{remark}

\begin{remark}
  The uniform coercivity condition (\ref{uniformcoercivity}) is a very
  reasonable condition although very difficult to prove in a framework
  as general as the one presented here. What is however clear in view
  of its definition, is that for all $\xi \in \Vs_0$ the matrix
  $\GG(\xi)$ is positive definite.
\end{remark}

We will not pursue further the analysis of problem~(\ref{eq:optimal-swimming}),
which requires analytical verification of conditions such as~(\ref{controllabilityhyp})
and~(\ref{uniformcoercivity}) on a case-by-case basis.
In the rest of our work, we will focus instead on
strategies to find solutions to  the optimization
problem~(\ref{eq:optimal-swimming}) numerically.

\section{Numerical Approximation}
\label{sec:bound-integr-appr}

The two main ingredients of the computer simulation of the swimmer
problem are an efficient numerical implementation of the Dirichlet
to Neumann map~\eqref{eq:dirichlet-to-neumann-map} and a strategy for
the search of minima of $\EE(\xi)$ in
Eq.~(\ref{eq:optimal-swimming}) satisfying the constraints.

We solve the first problem using the boundary element method (with a
custom code written in \texttt{C++}, using the \texttt{deal.II}
library~\cite{BangerthHartmannKanschat-2007-a,BangerthHartmannKanschat--a}
and based on the formulation presented
in~Ref.~\cite{Pozrikidis-1992-a}). We address the solution to
the constrained minimization problem using the reduced space
successive quadratic programming strategy
(rSQP)\cite{NocedalWright-1999-a} provided with the \texttt{Moocho}
package of the \texttt{Trilinos} \texttt{C++}
library.\cite{HerouxBartlettHowle-2005-a}

\subsection{Axisymmetric Boundary Integral Equations for Stokes Flow}
\label{sec:axisymm-bound-integr}

Our boundary integral formulation for axisymmetric Stokes flow follows
closely~Ref.~\cite{Pozrikidis-1992-a}:  if we consider Stokes equation in
free-space associated with a forcing term due to a Dirac mass centered
in $\vv y$ and weighted with the force vector $\vv b$, i.e.,
\begin{equation}
  \label{eq:free-space-stokes}
  \begin{aligned}
    & - \eta \Delta \u(\x) + \nabla p(\x) = -\nabla \cdot \vv \sigma  = \vv b\dirac(\x - \vv y)
    \qquad && \mbox{ in } \Re^3\\
    & \nabla\cdot \u(\x) =0 && \mbox{ in } \Re^3\\
    & \lim_{|\x|\to\infty} |\u(\x)| = 0 \\
    & \lim_{|\x|\to\infty} |p(\x)| = 0,
  \end{aligned}
\end{equation}
we can express the solutions $\u$, $p$ and $\vv \sigma$ using the
free-space Green's functions $\tilde{\vv
  {\mathcal S}}$, $\tilde{\vv P}$ and $\tilde{\vv T}$:
\begin{equation}
  \begin{split}
    \u(\x) & = \frac{1}{8\pi \eta}\tilde{\vv {\mathcal S}}(\x - \y) \vv b \\
    p(\x) & =  \frac{1}{8\pi } \tilde{\vv P}(\x - \y) \cdot \vv b\\
    \vv \sigma(\x) & =  \frac{1}{8\pi } \tilde{\vv T}(\x - \y) \vv b,
\end{split}
\label{eq:velocity-stokeslet}
\end{equation}
where
\begin{equation}
  \label{eq:free-space-green-functions}
  \begin{aligned}
    \tilde {\mathcal S}_{ij}(\vv r) & =  \frac{\delta_{ij}}{|\vv r|} +
    \frac{r_i r_j}{|\vv r|^3}\\
    \tilde P_i(\vv r) & = 2\frac{r_i}{|\vv r|^3}\\
    \tilde T_{ijk}(\vv r) & = -6 \frac{r_ir_jr_k}{|\vv r|^5}.
  \end{aligned}
\end{equation}

The axisymmetric approximation allows us to reduce the swimmer
problem to a boundary integral equation on a one-dimensional
curve, i.e., the swimmer is allowed to generate only
axial-symmetric surface forces, which are then transmitted to the
fluid.

We introduce the following notation for integrals on the surface
$\Gamma$, identified by the configuration $\X(s)=(X(s),\sigma(s)), \,\,s\in[0,1]$:
\begin{equation}
  \label{eq:change-of-coordinates}
  \begin{split}
    \int_{\Gamma} \tilde g \d \Gamma
    = & \ 2\pi \int_0^1 g(s) J_X(s) \sigma(s) \d s \\
    =: & \int_{\X(s)} g(s) \d s,
  \end{split}
\end{equation}
where we indicated with $\displaystyle J_X(s) = \left|\frac{\partial X}{\partial
    s}(s)\right|$, and the last line serves only to ease the
notation.

To solve Stokes equations in the axisymmetric case, we integrate
analytically the singular integral
Eq.~(\ref{eq:integral-equation}) around the axis of rotation,
while keeping the point $\x$ on the half-plane $\theta = 0$.

The azimuthal integration gives rise to a single layer potential
kernel, which can be expressed using complete elliptic integrals of
the first and second kind and produces a linearly coupled set of
one-dimensional integral equations.

We can write the resulting equation as
\begin{equation}
  \label{eq:integral-equation-symmetry}
  8\pi\eta  \v_\alpha(q) =
  - \int_{\X(s)} \vv{\mathcal S}_{\alpha \beta}(\X(q),\X(s)) \f_\beta(s) \d s,
\end{equation}
where the indices $\alpha$ and $\beta$ are either $x$ or $\sigma$, to
indicate the distance along the axis of symmetry, or the distance from
the axis itself.

The axisymmetric
kernel $\vv {\mathcal S}$ represents the single layer integrals of a \emph{ring}
of singularities passing through the point
$\tilde \X(s,0) = \Rot(0) \X(s) $:
\begin{equation}
  \label{eq:axisymmetric-stokes-kernel}
  \vv {\mathcal S} (\X,\Y):=  \int_0^{2\pi}
  \left[
    \begin{array}{cc}
      \tilde {\mathcal S}_{xx}\argS & (\tilde {\mathcal S}_{yx} \argS\cos\theta + \tilde {\mathcal S}_{zx}\argS \sin\theta) \\
      \tilde {\mathcal S}_{xy}\argS & (\tilde {\mathcal S}_{yy} \argS\cos\theta + \tilde {\mathcal S}_{zy}\argS \sin\theta)
    \end{array}
  \right] \d \theta,
\end{equation}
where $\tilde{\vv {\mathcal S}}\argS = \tilde{\vv {\mathcal
    S}}(\Rot(0) \X, \Rot(\theta) \Y)$ is the three-dimensional
\emph{Stokeslet} given in~(\ref{eq:free-space-green-functions}),
evaluated at the points $\Rot(0) \X$ and $\Rot(\theta) \Y$.

The above integration can be expressed in terms of complete elliptic
integrals of the first and second kind, which are defined,
respectively, as
\begin{equation}
  \label{eq:complete_elliptic_first}
  K(x) := \int_0^{\pi/2}{\frac{\d \theta}{\sqrt{1-x^2\sin^2\theta}}\d \alpha}, \qquad x \in [0,1),
\end{equation}
which is a logarithmically singular integral, and
\begin{equation}
  \label{eq:complete_elliptic_second}
  E(x) := \int_0^{\pi/2} \sqrt{1-x^2\sin^2\theta} \d \theta, \qquad x \in [0,1],
\end{equation}
which is a bounded integral.

Indeed, if we define
\begin{equation*}
  \X(q) = \left(\begin{array}{c}
      x_0 \\
      \sigma_0
    \end{array}
  \right),
  \qquad
  r := \sqrt{(x-x_0)^2+(\sigma-\sigma_0)^2},
  \qquad
  \Delta x = x-x_0,
\end{equation*}
and
\begin{equation*}
  k^2 = \frac{4\sigma \sigma_0}{(x-x_0)^2 + (\sigma +\sigma_0)^2},
\end{equation*}
then the axisymmetric Stokeslet is given by (see~Ref.~\cite{Pozrikidis-1992-a})
\begin{equation}
  \label{eq:axy-stokeslet-elliptic}
  \begin{split}
    {\mathcal S}_{xx} & = 2k \sqrt{\frac{\sigma}{\sigma_0}} \left( K(k) + \Delta x^2 \frac{E(k)}{r^2}\right) \\
    {\mathcal S}_{x\sigma} & = +\frac{k\Delta x}{\sigma}\sqrt{\frac{\sigma}{\sigma_0}} \left(K(k)-(\sigma_0^2-\sigma^2+\Delta x^2)\frac{E(k)}{r^2}\right) \\
    {\mathcal S}_{\sigma x} & = -\frac{k \Delta x}{\sigma_0} \sqrt{\frac{\sigma}{\sigma_0}}\left(K(k)+(\sigma_0^2-\sigma^2-\Delta x^2)\frac{E(k)}{r^2}\right) \\
    {\mathcal S}_{\sigma \sigma} & = \frac{k}{\sigma_0\sigma}\sqrt{\frac{\sigma}{\sigma_0}} \Bigg(
    (\sigma_0^2+\sigma^2+2\Delta x^2) K(k) - \\
    & \phantom{\frac{k}{\sigma_0\sigma}\sqrt{\frac{\sigma}{\sigma_0}} \Bigg( } \Big(2\Delta x^4
      +3\Delta x^2(\sigma_0^2+\sigma^2) + (\sigma_0^2-\sigma^2)^2\Big)\frac{E(k)}{r^2}\Bigg),
  \end{split}
\end{equation}
where the diagonal entries are singular and behave like ${\mathcal
  S}_{xx} \sim {\mathcal S}_{\sigma\sigma} \sim -2 \ln(r)$.

\subsection{Spatial Discretization}
\label{sec:discretization}

We employ an iso-parametric finite dimensional representation for the
spatial discretization of the configuration $\X$, and of the velocity
and force density $\u$ and $\f$. We call
$V_h=\text{span}\{\sbsp^i\}_{1\leq i \leq M}$ the finite dimensional space that we
use for the discretization. To ensure regularity and smoothness of the
boundary $\Gamma$ for any $t$, we choose a discretization based on
cubic bell splines.

In particular at any given time we identify elements of $V_h$ with the
coefficients of $M$ dimensional vectors, as in
\begin{equation}
  \label{eq:definitiona-finite-dimensional-vectors}
  \begin{split}
    \X_h(s) & = \X^i\sbsp^i(s)\\
    \v_h(s) & = \v^i\sbsp^i(s),
  \end{split}
\end{equation}
where we used the summation convention, and $i$ goes from 1 to $M$.
The superscripts here indicate the indices of the basis functions that
span the space $V_h$. When confusion cannot arise, we will use $\X_h$
without ``$(s)$'' to indicate the $\Re^M$ vector of components $\X^i,\,
i=1\dots M$ which uniquely identify the vector valued function
$\X_h(s)$ in $V_h$ as
in~(\ref{eq:definitiona-finite-dimensional-vectors}).

To take into account collections of shapes, it is possible to
introduce controlled discontinuities in the basis functions
$\sbsp^i(s)$, by overlapping a sufficient number of nodes in the
knot-span that generates the B-Spline basis, as explained, for
example, in Ref.~\cite{PieglTiller-1997-a}.

The Galerkin approach to the boundary element method can be understood
by multiplying Eq.~(\ref{eq:integral-equation-symmetry}) on both
sides with test functions (or virtual forces) $\sbsp$ and integrating
on $\Gamma$
\begin{multline}
  \label{eq:integral-equation-symmetry-variational}
  \int_{\X(s)} \v(s) \cdot \sbsp(s) \d s = \\
  - \int_{\X(q)} \sbsp(q) \cdot \frac{1}{8\pi\eta} \int_{\X(s)} \vv {\mathcal
    S}(\X(q),\X(s)) \f(s) \d s \d q \\
  \forall \sbsp \in C^\infty([0,1])^2,
\end{multline}
which, on a discrete level, yields the discrete
version of the Dirichlet to Neumann map:
\begin{equation}
  \label{eq:discrete-DN-map}
  \f_h  = \A^{-1} \M \u_h,
\end{equation}
where both matrices $\A$ and $\M$ depend on $\X_h(s)$ and are defined as
\begin{equation}
  \label{eq:discrete-A}
  \A^{ij} = - \int_{\X_h(q)} \sbsp^i(q) \cdot \frac{1}{8\pi\eta} \int_{\X_h(s)} \vv {\mathcal
    S}(\X_h(q),\X_h(s)) \sbsp^j(s) \d s \d q,
\end{equation}
and
\begin{equation}
  \label{eq:discrete-M}
  \M^{ij} = \int_{\X_h(s)} \sbsp^j(s) \cdot \sbsp^i(s) \d s.
\end{equation}

Numerical computation of the matrix entries of both $\A$ and $\M$ is
performed using high order Gauss quadrature formulae, as well as
Gauss-Log quadrature formulae for the diagonal entries of $\A$, in
order to properly take into account the logarithmic singularity of the
axisymmetric Stokeslet (see, for example,
Ref.~\cite{Pozrikidis-1992-a,Pozrikidis-2002-a,Roumeliotis-2000-a}).

Notice that, once we have two vectors $\f_h$ and $\v_h$ in $\Re^M$
that represent objects of $V_h$, we can compute their integrals on
$\X_h(s)$ by scalar product in $\Re^M$ through the matrix $\M$:
\begin{equation}
  \label{eq:def-integral-Re-M}
  \int_{\X(s)} \f_h(s)\cdot \v_h(s) \d s= \f_h \cdot \M\v_h = \v_h
  \cdot \M\f_h = \v^i \M^{ij} \f^j,
\end{equation}
where in the left hand side $(\cdot)$ has the meaning of scalar
product between vector \emph{functions} of $\Re^2$, while on the right
hand sides $(\cdot)$ should be read as the scalar product in $\Re^M$.

Once we know how to discretize a domain $\Gamma$ and write a discrete
version of the Dirichlet to Neumann Map \emph{for any configuration},
we simply define the \emph{discrete shape basis functions} as the
counter part of $\partial \bar \X(\xi,\varphi)/\partial\xi$, i.e.,
 \begin{equation}
  \label{eq:discrete-shape-basis}
  \begin{split}
    & \basexi_h^i(s) = \frac{\partial \bar \X_h(\xi, s)}{\partial \xi_i}
    \qquad i=1\dots N\\
    & \axis_h(s) = \axis(s),
  \end{split}
\end{equation}
where $\bar \X(\xi, s)$ indicates a point on the $(x,\sigma)-$half-plane
 of the configuration $\Gamma_0(\xi)$.

The discrete shape basis functions $\basexi_h^i(s)$ are simply
constructed as the b-spline interpolation of their continuous counter
parts, which are a datum of the problem, while the discrete basis
function $\axis_h(s)$ for the displacement vector space is simply
equal to the discrete representation of the constant unit vector
directed along the axis of symmetry.

Given $N+1$ arbitrary time functions $(\xi(t), \varphi(t))$, it is
then possible to build the discrete configuration $\X_h(\xi(t),
\varphi(t))$, the discrete basis functions for the rate of shape
change $\basexi^i_h$ and to compute the Dirichlet to Neumann map for
arbitrary data defined on the discrete curve.

In particular we can compute the discrete versions of $\VV$ and $\GG$
as
\begin{equation}
  \label{eq:def-NNh-MMh}
   \VV_{i,h} (\xi)  :=  - \frac{ (\A^{-1} \M \basexi_h^i(\xi) )\cdot
      \M\axis_h } { (\A^{-1} \M \axis_h) \cdot \M\axis_h },
\end{equation}
and
\begin{equation}
  \label{eq:def-GGh}
   \GG_{ij,h} (\xi)  :=   \big(\A^{-1} \M ( \basexi_h^i(\xi) +
   \VV_{i,h}(\xi)\axis_h) \big)\cdot \M( \basexi_h^j(\xi) +
   \VV_{j,h}(\xi)\axis_h).
\end{equation}

The semi-discrete formulation of the optimal swimming problem is
simply obtained by substituting in
Eq.~(\ref{eq:optimal-swimming}) $\VV$ and $\GG$ with their
\emph{spatial} discrete counterparts:
\begin{subequations}
  \label{eq:optimal-swimming-semi-discrete}
  \begin{align}
    \min_{\xi \in \Vs}\qquad & \int_0^T \GG_h(\xi)\dot \xi \cdot \dot
    \xi \d t =: \EE_h(\xi) \label{eq:opt-minimality-semi-discrete}\\
    \text{subject to }\qquad &\int_0^T\VV_h(\xi)\cdot\dot \xi \d t- c
    = : \CC_h^c(\xi) = 0 \label{eq:opt-feasibility-semi-discrete}.
  \end{align}
\end{subequations}

\subsection{Time Discretization}
\label{sec:time-discretization}

The dependency of the domain $\Gamma$ (or, equivalently, of its
axisymmetric representation $\X$) on the time variable $t$ is
\emph{only} through the \emph{shape} and \emph{position} functions
$\xi_i(t)$ and $\varphi(t)$. In particular we rewrote the problem in a
way that makes it only dependent on the shape parameters $\xi \in
\Vs$.

We approximate also these functions with cubic b-spline functions,
with a possibly different dimension $Q$. In particular we have
\begin{equation}
  \label{eq:definition-Vsh}
  \Vs_h= \left\{ \xi_h(t) \in
  \left(\text{span}\{\tbsp^\alpha(t)\}_{\alpha=1}^Q\right)^N,
  \text{ s.t. } \xi_h(0) = \xi_h(T) \right\},
\end{equation}
where we use Greek letters to
distinguish basis functions referring to the time variable $t$ from
those referring to the space variable $s$ and again we can write
\begin{equation}
  \label{eq:definition-finite-dimensionale-xi-phi}
  \begin{split}
    &\xi_h(t)_i = \xi_i^\alpha \tbsp^\alpha(t) \qquad i = 1\ldots N,
  \end{split}
\end{equation}
where summation is implied on the repeated indices $\alpha$ which go
from $1$ to $Q$.

Notice that a numerical approximation $\xi_h$ is now a vector of
$\Re^{N\times Q}$, which we identify either with two indices
$\xi_i^\alpha$, or with a single capital letter index $I=1\dots N\times
Q$, where
\begin{equation}
  \label{eq:indices-reordering}
  I := i + (\alpha-1)N, \qquad \xi_h :=
  \left\{\{\xi_i^\alpha\}_{i=1}^{N} \right\}_{\alpha=1}^Q = \{\xi_I\}_{I=1}^{N\times Q},
\end{equation}
that is
\begin{equation}
  \label{eq:xih-unrolled}
  \xi_h := \left[
    \xi^1_1, \xi^1_2,\dots,\xi^1_N, \quad
    \xi^2_1, \xi^2_2, \dots,\xi^2_N, \quad \dots, \quad
    \xi^Q_1, \xi^Q_2,\dots,\xi^Q_N\right]'.
\end{equation}

\subsection{Discrete Optimal Swimming}
\label{sec:discr-optim-swimm}

With the discretization $\Vs_h$ of the infinite dimensional space
$\Vs$, we can reduce the minimization problem to a finite dimensional
one defined on $\Re^{Q\times N}$, where $Q$ is the total number of
basis splines selected for the discretization of the $N$ shape
parameters in time, and the fully discrete problem becomes
\begin{subequations}
  \label{eq:optimal-swimming-fully-discrete}
  \begin{align}
    \min_{\xi_h \in \Vs_h}\qquad &   \EE_h(\xi_h)  :=  \Eh(\xi_h) \xi_h \cdot \xi_h
    \label{eq:opt-minimality-fully-discrete}\\
    \text{subject to }\qquad & \CC_h^c(\xi_h) := \Ch(\xi_h) \cdot \xi_h -c
    = 0 \label{eq:opt-feasibility-fully-discrete},
  \end{align}
\end{subequations}
where the matrix $\Eh(\xi_h)$ in $\Re^{(N\times Q)\times (N\times
  Q)}$ and the vector $\Ch(\xi_h)$ in $\Re^{(N\times Q)}$ are
  defined as
\begin{equation}
  \label{eq:definition-minimality-feasability-matrices}
  \begin{split}
    \Eh_{IJ}(\xi_h) := & \int_0^T \GG_{ij}(\xi_h(t) ) \dot
    \tbsp^\alpha(t) \cdot \dot \tbsp^\beta(t) \d t\\
    \Ch_{I}(\xi_h) := & \int_0^T \VV_i(\xi_h(t) ) \dot \tbsp^\beta(t) \d t,
  \end{split}
\end{equation}
and the relationship between the indices $I,J$ and $i,j,\alpha, \beta$
is given by (\ref{eq:indices-reordering}) and $J = j+(\beta-1)N$.

The problem is now written in the framework of classical finite
dimensional optimization. Indeed, if we define the nabla operator
$\nabla$ applied to a function of the vector variable $\xi_h$ as
\begin{equation}
  \label{eq:definition-nabla}
  (\nabla \vv F(\xi_h))_I := \frac{\partial \vv F(\xi_h)}{\partial \xi_I},
\end{equation}
we can then write the Lagrangian $\LL(\xi_h, \lambda)$ and its
derivatives associated with
problem~(\ref{eq:optimal-swimming-fully-discrete}) as
\begin{equation}
  \label{eq:lagrangian}
  \begin{split}
    \LL(\xi_h, \lambda) := & \EE_h(\xi_h) + \lambda
    \CC_h^c(\xi_h) \\
     \nabla \LL(\xi_h, \lambda) := &  \nabla \EE_h(\xi_h)+
     \lambda \nabla\CC_h^c(\xi_h)\\
     \nabla^2 \LL(\xi_h, \lambda) := & \nabla^2\EE_h(\xi_h)+
     \lambda\nabla^2\CC_h^c(\xi_h).
  \end{split}
\end{equation}

The first and second-order necessary Karush-Kuhn-Tucker (KKT)
optimality conditions for a solution $(\bar \xi_h, \bar\lambda)$ to
(\ref{eq:optimal-swimming-fully-discrete}) are: \begin{equation}
  \label{eq:minimality-KKT-conditions-fully-discrete}
  \begin{split}
    &\nabla \LL(\bar\xi_h, \bar\lambda) =  \nabla \EE_h(\bar\xi_h) + \bar\lambda
    \nabla \CC_h^c(\bar\xi_h) = 0\\
    &\CC_h^c(\bar\xi_h) =  0 \\
    &\nabla^2 \LL(\bar\xi_h, \bar\lambda) \eta \cdot \eta \quad \ge  0
    \qquad \forall \eta \in \Vs_h,
  \end{split}
\end{equation}
also known as \emph{linear dependence of gradients},
\emph{feasibility} and \emph{curvature} conditions.

A popular class of methods for solving non linear constrained
minimization problems is successive quadratic programming (SQP) (see,
e.g., Ref.~\cite{NocedalWright-1999-a}), which is equivalent to applying
Newton's method to solve the the minimality
conditions~(\ref{eq:minimality-KKT-conditions-fully-discrete}).

At each Newton iteration $k$
for~(\ref{eq:minimality-KKT-conditions-fully-discrete}), the linear
subproblem (also known as the KKT system) takes the form
\begin{equation}
  \label{eq:KKT-system}
  \left(
    \begin{array}{cc}
      W & A\\
      A^T &  0
    \end{array}
    \right)
    \left(
      \begin{array}{c}
        d\\
        d_\lambda
      \end{array}
    \right) = -
    \left(
    \begin{array}{c}
      \nabla \LL^{(k)}\\
      \CC_c^{(k)}
    \end{array}
    \right),
\end{equation}
where
\begin{equation}
  \label{eq:minimality-KKT-conditions-h}
  \begin{aligned}
    d &:= \xi_h^{(k+1)}-\xi_h^{(k)} \\
    d_\lambda & := \lambda^{(k+1)}-\lambda^{(k)}\\
    W &:= \nabla^2\LL(\xi_h^{(k)},  \lambda^{(k)})\\
    A &:= \nabla \CC_h^c(\xi_h^{(k)}) \\
    \LL^{(k)} & := \LL(\xi_h^{(k)},  \lambda^{(k)}) \\
    \CC_c^{(k)} & := \CC^c_h(\xi_h^{(k)} ),
  \end{aligned}
\end{equation}
and we use superscripts between parenthesis $(k)$ to refer to the
$k$-th Newton iteration, not to be confused with the index of the
b-spline basis functions.

Once we obtained a new estimate of the solution $(\xi_h^{(k+1)},
\lambda^{(k+1)})$, the error in the optimality
conditions~(\ref{eq:minimality-KKT-conditions-fully-discrete}) is
checked. If these KKT errors are within some specified tolerance, the
algorithm is terminated with the optimal solution.

If the KKT error is too large, the functions and gradients are then
computed at the new point $\xi_h^{(k+1)}$ and another KKT
subproblem~(\ref{eq:KKT-system}) is solved which generates another
increment $d$, until convergence or failure.

A successful application of the SQP method requires one to provide
exact information about the hessian of the Lagrangian $W$. While this
ensures second order convergence of the SQP method, it is not feasible
in many practical cases, including ours, where the direct computation
of $W$ is overly expensive.

To address this issue, we use a reduced-space SQP (rSQP) method. In
rSQP methods, the full-space QP subproblem~(\ref{eq:KKT-system}) is
decomposed into two smaller subproblems, where the optimization vector
variable $\xi_h$ is split into \emph{dependent} and \emph{independent}
variables, by using the linearized equality constraint as a mean to
express the dependent variables (often referred to as the
\emph{control} variables) in terms of the independent ones (also
referred to as the \emph{state} variables).

Using this decomposition, the two subproblems can be solved
effectively by using an approximation of $W$. We use the implementation
of the rSQP method included in the \texttt{Moocho} package of the
\texttt{Trilinos} library.\cite{Bartlett-2009-a}

\section{Numerical Results}
\label{sec:numerical-results}

The goal of this section is to present some results obtained with
our code. In section~\ref{sec:three-sphere-swimmer} we present our
results for the optimal stroke of the Golestanian swimmer. These
should be confronted with the results already presented in
Ref.~\cite{AlougesDeSimoneLefebvre-2007-a} and
Ref.~\cite{AlougesDeSimoneLefebvre-2008-a} with a completely different
approach both for the resolution of Stokes equations and for the
minimization of the expended energy, which we take as our
reference results.

While from the quantitative point of view we observe some marginal
differences, qualitatively the method presented in this paper is
capable of reproducing the optimal strokes obtained in the previously
mentioned works.

We should mention here that the main differences lay in the fact that
in those works the approximation of the $\VV(\xi)$ and $\GG (\xi)$
fields was performed using a Finite Element Method for the simulation
of axisymmetric Stokes Flow in a box that was taken large compared to
the dimension of the swimmer, but of \emph{finite} dimension. In
comparison, the method we present here is based on the use of Boundary
Integral Approximation, which should yield more accurate results on
infinite domains such as the one we are interested on.

As a second remark, we observe that the optimal stroke in
Ref.~\cite{AlougesDeSimoneLefebvre-2007-a} and
Ref.~\cite{AlougesDeSimoneLefebvre-2008-a} is obtained by writing
explicitly the Euler-Lagrange equation for the constrained
optimization problem, and by solving the resulting system of ODEs with
an explicit Runge Kutta method, coupled with a shooting approach in
order to enforce the periodicity of the strokes.

In this paper, on the other hand, the solution of the minimization
problem itself is left to an external library which explores the space
of shapes in a very efficient and flexible way, allowing us to solve
problems that could not be addressed before. In particular we can
solve problems in which we include inequality constraints on the shape
variables $\xi$, and we can release the constraints on the initial and
final shape, letting the minimizer find the optimal starting shape for
our swimmers.

In the experiments that follow, we use water at room temperature
($20^{o}$ C) as the surrounding fluid, and we express lengths in
millimeters ($mm$), time in seconds ($s$) and weight in milligrams
($mg$). Using these units, the viscosity of water is approximately one
($1 mPa\, s = 1 mg\, mm^{-1}\, s^{-1}$), and the energy is expressed in
pico-Joules ($1pJ=1mg\, mm^2\, s^{-2}=10^{-12} J$).

The tests were performed on a Macbook Pro, with 2.16 GHz Intel Core 2
Duo processor and 2 GB 667 MHz DDR2 SDRAM. The average running time
for the experiments presented in Sec.~\ref{sec:three-sphere-swimmer}
is about 30 seconds, with a maximum memory consumption of 64MB, while
the experiments of Sec.~\ref{sec:stick-donut-swimmer} required on
average 10 seconds, with roughly the same memory consumption.

The low running time is possible thanks to the use of cubic
B-Spline, which allows one to obtain accurate solutions using very few
degrees of freedom. In all the experiments we present, we used 15
control points on each portion of the swimmer for the spatial
discretization (i.e., 45 total control points for the three sphere
swimmer and 30 for the stick and donut swimmer) and 15 control points
for the time discretization of each shape and position variable. 

\subsection{Three Sphere Swimmer}
\label{sec:three-sphere-swimmer}

The three sphere swimmer is among the simplest axisymmetric swimmers
one can think of. It consists of three linked spheres that can only
vary their reciprocal distance. The swimmer is completely defined once
we have the radius of the three spheres (which we assume to be the
same for simplicity), and the location of the three centers of the
spheres on the axis of the movement (one size parameter and three
positional variables).

\begin{figure}
  \centering
  \input{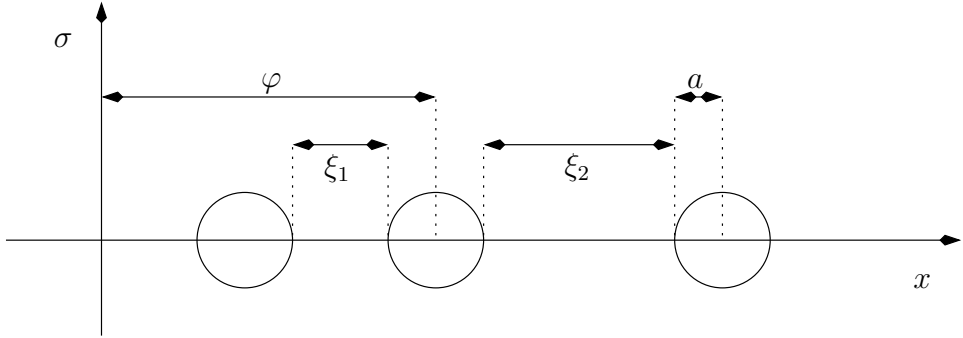}
  \caption{Definition of the Three Sphere Swimmer.}
  \label{fig:three-sphere-swimmer-def}
\end{figure}

Equivalently, we can use the more convenient representation given by
the location of the center of the central sphere ($\varphi$) and the
distances between the two lateral spheres and the central one ($\xi_1$
and $\xi_2$), as shown in
Figure~\ref{fig:three-sphere-swimmer-def}. These two different
representations are clearly equivalent (see, e.g.,
Ref.~\cite{AlougesDeSimoneLefebvre-2007-a}).

Using this convention $\X(\xi, \varphi, s)$ can be parametrized as
\begin{equation}
  \label{eq:definition-X-3spheres}
  \begin{aligned}
    \X(\xi, \varphi, s) = & \left(
      \begin{array}{c}
        \varphi-\xi_1 -a(2+\cos(3\pi s)) \\
        a\sin(3 \pi s)
      \end{array}
    \right)  && \text{ if } s \in [0, 1/3) \\
    & \\
   \X(\xi, \varphi, s) = & \left(
      \begin{array}{c}
        \varphi -a\cos(3\pi (s-1/3)) \\
        a\sin(3 \pi (s-1/3))
      \end{array}
    \right)  && \text{ if } s \in [1/3, 2/3] \\
    & \\
  \X(\xi, \varphi, s) = & \left(
      \begin{array}{c}
        \varphi+\xi_2 + a(2-\cos(3\pi s)) \\
        a\sin(3\pi (s-2/3))
      \end{array}
    \right)  && \text{ if } s \in (2/3,1],
  \end{aligned} \end{equation}
which implies that the basis functions for shape changes $\basexi_i :=
\partial \X/\partial \xi_i$ are defined as
\begin{equation}
  \label{eq:basis-xi-3-spheres}
  \basexi_1(s) = -\chi_{[0, 1/3)}(s)\axis, \qquad   \basexi_2(s) = \chi_{(2/3, 1]}(s)\axis,
\end{equation}
where $\chi_{A}(s)$ is the function which is one if $s\in A$, and zero
otherwise. We constrain the shape variables $\xi_i$ to be in the
interval $[0, 6a]$. Notice that in this paper $\xi$ is the touching
distance between the spheres, and not the distance between the centers
of the spheres, as in~Ref.~\cite{AlougesDeSimoneLefebvre-2007-a}.

\begin{figure}[!htb]
  \centering
  \includegraphics[width=\textwidth]{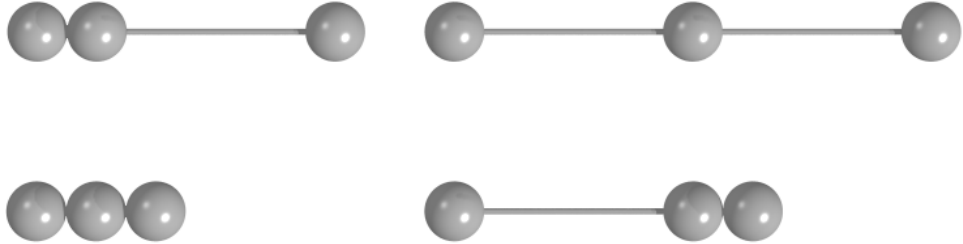}
  \caption{Extremal shape configurations for a Three-Sphere swimmer of
    radius $a=.05mm$,  from left to right and from bottom to top:
    $\xi = (0,0)$, $\xi = (6a,0)$, $\xi=(0,6a)$ and $\xi=(6a,6a)$.}
  \label{fig:extremal-shapes-3-spheres}
\end{figure}

In Figure~\ref{fig:extremal-shapes-3-spheres} we show the corners of
the box $[0, 6a]^2$ in which the three sphere swimmer shape is
constrained to stay, for a swimmer of radius $a=.05mm$.

The basis function for the change in position, $\partial \X/\partial
\varphi$ is simply equal to the unit vector $\axis$. In
Figure~\ref{fig:basis-u-e}, we show the three basis functions that
allow one to fully describe the system, together with their Dirichlet
to Neumann maps, evaluated at the configuration $\xi = (.05mm, .05mm)$,
which is equivalent to asking that the distance between the centers of
the spheres is $3a$.

\begin{figure}[!htb]
  \centering
  \includegraphics[width=.45\textwidth]{\Figpath{3sphere-u0}}
  \hfill
  \includegraphics[width=.45\textwidth]{\Figpath{3sphere-f0}}

  \includegraphics[width=.45\textwidth]{\Figpath{3sphere-u1}}
  \hfill
  \includegraphics[width=.45\textwidth]{\Figpath{3sphere-f1}}

  \includegraphics[width=.45\textwidth]{\Figpath{3sphere-phi-dot}}
  \hfill
  \includegraphics[width=.45\textwidth]{\Figpath{3sphere-f-phi-dot}}

  \caption{Basis function $\basexi_i(s)$ and $\axis$ (left) and their
    Dirichlet to Neumann  maps (right) for $\xi = (a, a)$. }
  \label{fig:basis-u-e}
\end{figure}

\begin{figure}[!htb]
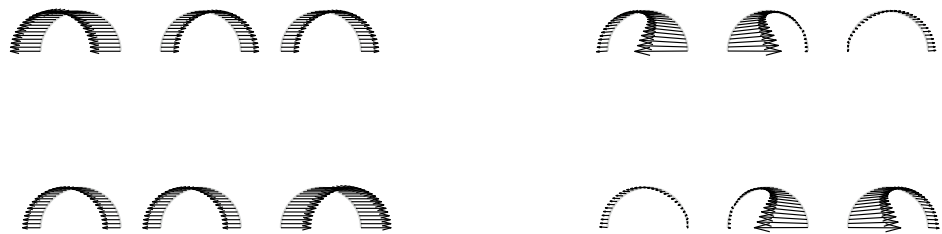

  \centering
  \includegraphics[width=.45\textwidth]{\Figpath{3sphere-w0}}
  \hfill
  \includegraphics[width=.45\textwidth]{\Figpath{3sphere-f-w0}}

  \includegraphics[width=.45\textwidth]{\Figpath{3sphere-w1}}
  \hfill
  \includegraphics[width=.45\textwidth]{\Figpath{3sphere-f-w1}}

  \caption{The force-free basis function $\basev_i(s)$ (left) and their
    Dirichlet to Neumann maps (right) for $\xi = (a, a)$. }
  \label{fig:basis-w}
\end{figure}

Once we have the various basis functions, it is easy to combine them
linearly and obtain a basis for \emph{force-free} movement. This
is what Figure~\ref{fig:basis-w} shows, where the basis functions
$\basev_i = \basexi_i+\VV_i\axis$ are plotted, together with their
Dirichlet to Neumann map. Notice that, by construction, the
integral of the forces on the configuration $\Gamma$ yields zero
to machine precision.

Minimizing the expended energy and forcing the displacement to be
a given datum, we obtain the same qualitative results as
in~Ref.~\cite{AlougesDeSimoneLefebvre-2007-a}
and~Ref.~\cite{AlougesDeSimoneLefebvre-2008-a}. The approach we propose
in this paper, however, allows us to study also the problem of
finding the optimal stroke \emph{without} assigning an initial
(and final) shape.

\begin{figure}[!htb]
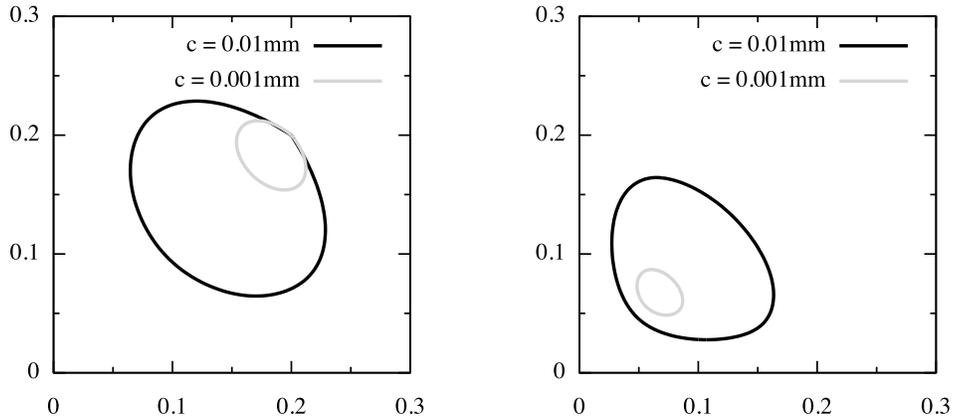

  \centering
  \includegraphics[width=.49\textwidth]{\Figpath{3sphere-shape}}
  \hfill
  \includegraphics[width=.49\textwidth]{\Figpath{3sphere-shape-no-initial}}
  \caption{Path in shape space describing optimal strokes for three sphere swimmer
    of radius $.05mm$, swimming  $.01mm$ and $.001mm$ in 1
    second, imposing the initial shape $\xi=(.3mm, .3mm)$ on the
    left, and without imposing an initial shape on the right.}
  \label{fig:3sphere-optimal-stroke}
\end{figure}

The closed paths in shape space for  target displacements of
$.01mm$ and $.001mm$ can be seen in
Figure~\ref{fig:3sphere-optimal-stroke} for both the case where
the initial shape is fixed to be $(.2mm, .2mm)$ and for the case
where no constraints are imposed.

\begin{figure}[!htb]
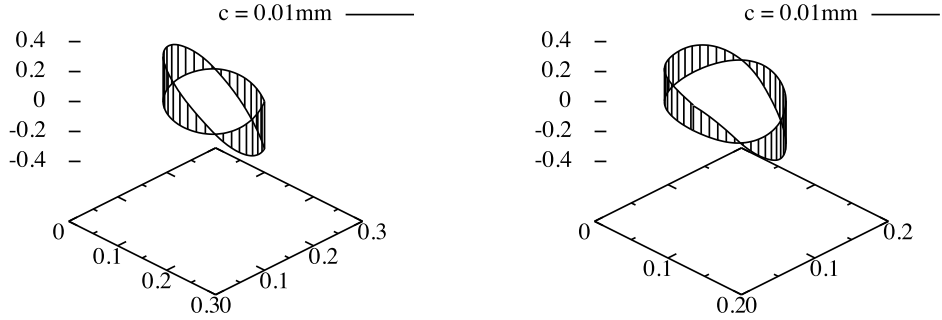

  \centering
  \includegraphics[width=.49\textwidth]{\Figpath{3sphere-phidot}}
  \hfill
  \includegraphics[width=.49\textwidth]{\Figpath{3sphere-phidot-no-initial}}
  \caption{Diagram of translational velocities for three sphere swimmer of radius
    $.05mm$ swimming $.01mm$ in 1 second imposing the initial
    shape $\xi=(.3mm, .3mm)$ on the left, and without imposing an initial
    shape on the right.}
  \label{fig:3sphere-optimal-phi-dot}
\end{figure}

In Figure~\ref{fig:3sphere-optimal-phi-dot} we can observe how the
optimal stroke can be separated into a propulsive part and a
recovery part, where the propulsive part is the one where the
velocity is bigger and positive, while the recovery part is the
one where the velocity is smaller and negative.

We would like to emphasize here that, in general, if we fix the
initial and final shapes, the velocity is free to be discontinuous at
$t=kT$ with $k$ integer, however the truly optimal stroke (the one
selected without impositions on the starting shape) is one where also
the velocity is continuous, as shown in
Figure~\ref{fig:3sphere-optimal-stroke} on the right.

\subsection{Stick and Donut swimmer}
\label{sec:stick-donut-swimmer}

We now present a new model swimmer, which simulates the swimming
mechanism of many biological organisms made of a body and a
propulsive apparatus consisting of appendages that change shape to
induce propulsion. A nice example is  \emph{Chlamydomonas
Reinhardtii}, a unicellular eukaryotic alga with a body size of
roughly $10$ $\mu$m,  that swims by executing a movement with its
two flexible flagella which is closely reminiscent of the breast
stroke of an olympic swimmer. An axisymmetric version of this
swimming style, though at larger scales, is that of the jelly
fish. In the model swimmer we propose, the body is schematized
with a cylinder capped with two half-spheres at the ends (the
\emph{Stick}) and the propulsive apparatus is schematized with a
torus with elliptical cross-section of variable major and minor
radii (the \emph{Donut}).

\begin{figure}[!htb]
  \centering
  \resizebox{\textwidth}{!}{
    \input{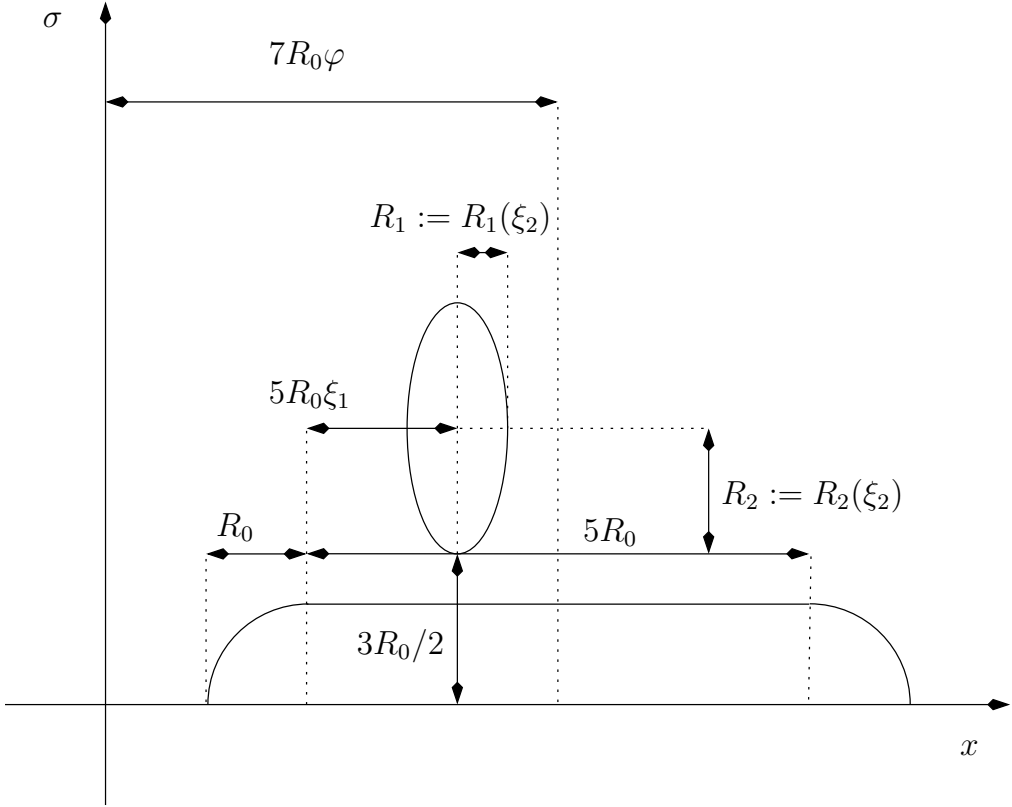}
  }
  \caption{Definition of the ``Stick and Donut'' swimmer.}
  \label{fig:stick-donut-def}
\end{figure}

Figure~\ref{fig:stick-donut-def} shows a section of the swimmer on
the $(x,\sigma)$-half plane, for $\theta = 0$. The radius of the
stick is set to be $R_0$, and all other dimensions are scaled with
respect to this one, so that the touching distance between the
stick and the donut is fixed to $3R_0/2$ and the volumes of the
stick and the donut are fixed and equal.

Swimming is achieved by translating the center of the donut along
the direction of the stick (non-dimensional shape variable
$\xi_1$, constrained in the interval $[0,1]$), so that the center
of the donut is always inside the ``body'' of the swimmer), and by
varying the radii of the donut section.

For the non-dimensional shape parameter $\xi_2=0$, the horizontal
radius $R_1$ of the donut section is equal to $5R_0/2$, i.e., it
has the same length of the stick without the half spheres, while
for $\xi_2=1$, the horizontal radius $R_1$ is equal to $R_0/5$.

The vertical radius of the swimmer and the vertical position of the
center of the ellipse section are adjusted automatically as a function
of $\xi_2$ in order to maintain the volume of the donut constantly
equal to the volume of the stick and to maintain the distance between
the donut and the stick constant.

The absolute position of the swimmer on the axis of symmetry $x$,
is given by the non-dimensional variable $\varphi$, so that when
$\varphi $ is equal to one, then the swimmer has moved of one
entire body length ($7R_0$).

\begin{figure}[!htb]
  \centering
  \includegraphics[width=.49\textwidth]{\Figpath{stick-donut-u0}}
  \hfill
  \includegraphics[width=.49\textwidth]{\Figpath{stick-donut-f0}}

  \includegraphics[width=.49\textwidth]{\Figpath{stick-donut-u1}}
  \hfill
  \includegraphics[width=.49\textwidth]{\Figpath{stick-donut-f1}}

  \includegraphics[width=.49\textwidth]{\Figpath{stick-donut-phi}}
  \hfill
  \includegraphics[width=.49\textwidth]{\Figpath{stick-donut-f-phi}}
  \caption{Basis functions $\basexi_i(s)$ and $\axis(s)$
    (left) and their Dirichlet to Neumann
    map $\f_1(s)$ and $\f_{\axis}$ (rigth) for $\xi = (.5, .5)$. }
  \label{fig:stick-donut-basis}
\end{figure}

\begin{figure}[!htb]
  \centering
  \includegraphics[width=.49\textwidth]{\Figpath{stick-donut-w0}}
  \hfill
  \includegraphics[width=.49\textwidth]{\Figpath{stick-donut-f-w0}}

  \includegraphics[width=.49\textwidth]{\Figpath{stick-donut-w1}}
  \hfill
  \includegraphics[width=.49\textwidth]{\Figpath{stick-donut-f-w1}}
 \caption{Force free basis functions $\basev_i(s)$
    (left) and their Dirichlet to Neumann
    map $\f_{\basev_i}$ (rigth) for $\xi = (.5, .5)$. }
  \label{fig:stick-donut-basis-force-free}
\end{figure}

Figures~\ref{fig:stick-donut-basis}
and~\ref{fig:stick-donut-basis-force-free} show the basis functions
$\basexi_i$ and $7R_0\axis$ and the corresponding force-free basis
functions $\basev_i$ with their Dirichlet to Neumann maps. From
Figure~\ref{fig:stick-donut-basis-force-free} it is evident how
changes in the shape parameter $\xi_2$ do not induce much displacement
along the axis of symmetry.

\begin{figure}[!htb]
  \centering
  \includegraphics[width=\textwidth]{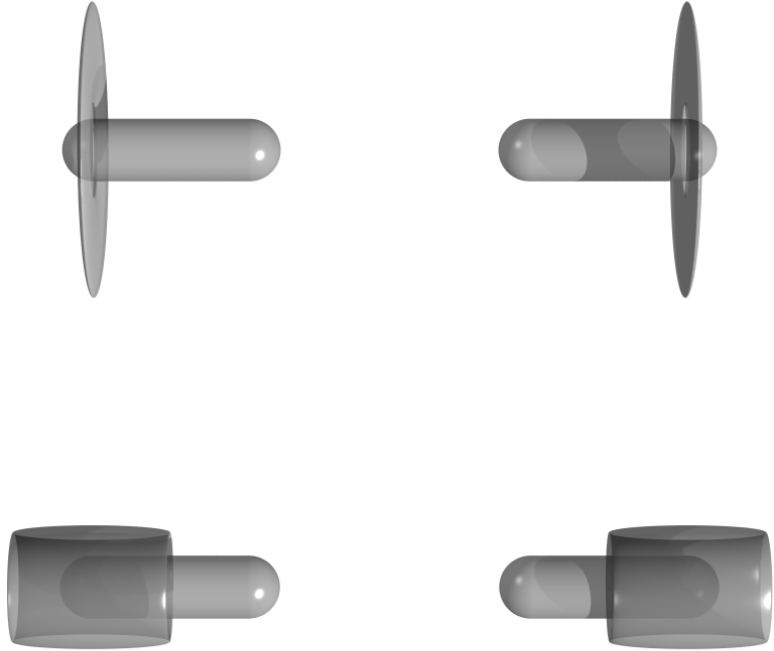}
  \caption{Extremal shape configurations for the stick and donut swimmer:
    from left to right and from bottom to top: $\xi = (0,0)$, $\xi =
    (1,0)$, $\xi=(0,1)$ and $\xi=(1,1)$.}
  \label{fig:extremal-shapes-stick-donut}
\end{figure}

\begin{figure}[!htb]
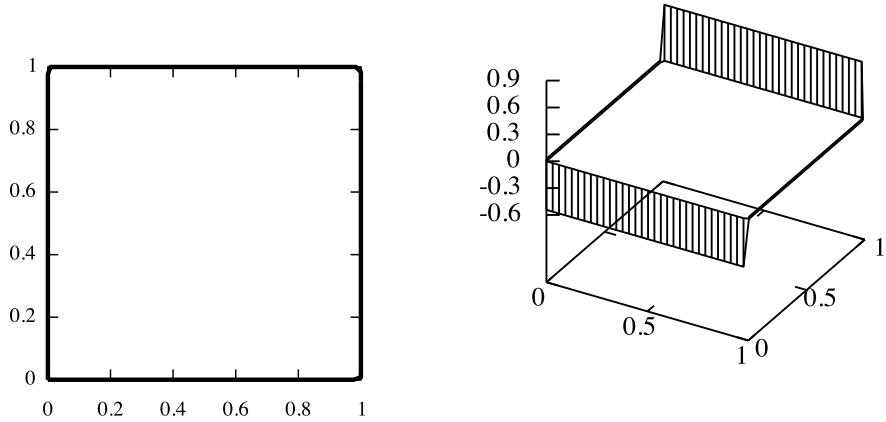

  \centering
  \includegraphics[width=.43\textwidth]{\Figpath{stick-shape-square}}
  \hfill
  \includegraphics[width=.55\textwidth]{\Figpath{stick-phi-dot-square}}
  \caption{Path in shape space describing a square strokes for stick
    and donut swimmer with base
    radius of $.034mm$, swimming $.034mm$ in 1
    second (left), with its propulsion diagram (right).}
  \label{fig:stick-square-stroke}
\end{figure}

We constrain the swimmer shape variables to be included in the
non-dimensional box $[0,1]^2$, whose corners are shown in
Figure~\ref{fig:extremal-shapes-stick-donut}.  A square stroke that
explores in a clockwise manner this shape space is presented in
Figure~\ref{fig:stick-square-stroke}. 

The fact that changes in the shape parameter $\xi_2$ have little
effect on the overall displacement of the swimmer can also be inferred
from the left and right sides of the square stroke, in the right part
of Figure~\ref{fig:stick-square-stroke}, which shows how the velocity
$\dot\varphi$ of the swimmer due to changes in $\xi_2$ are negligible
when compared with the velocity due to changes in $\xi_1$ (top and
bottom sides of the square stroke).

\begin{figure}[!htb]
  \centering
  \input{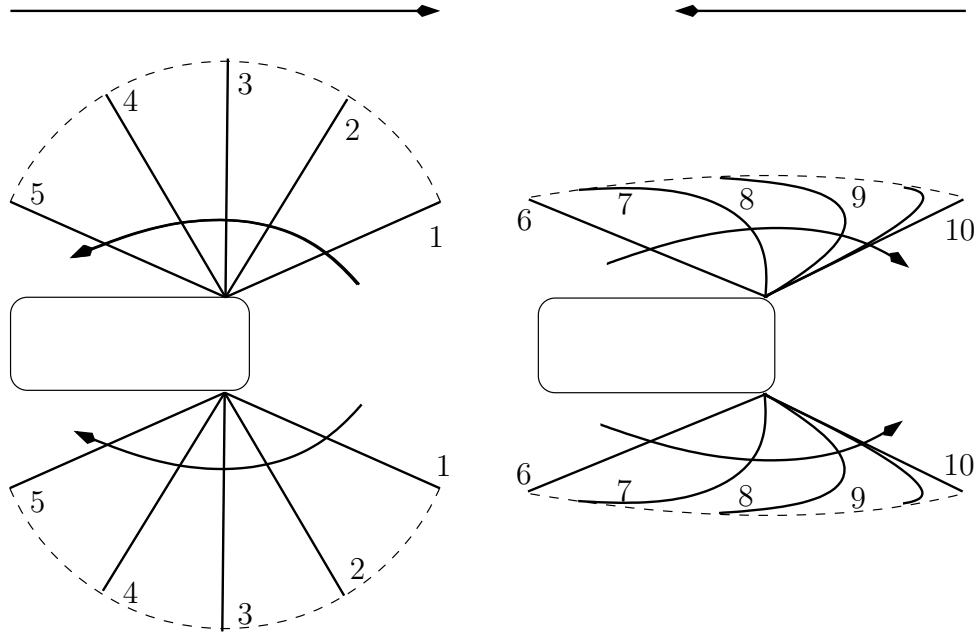}
  \caption{Ciliary stroke cycle. Power phase (left) and recovery
    phase (right).}
  \label{fig:ciliary-cicle}
\end{figure}

The separation of stroke cycles into a \emph{power} or
\emph{propulsive} phase, in which appendages have maximal
\emph{hydrodynamical resistance}, and a \emph{recovery} phase, in
which the swimmer tries to minimize viscous drag forces on its
propulsive appendages, are very common in nature. A typical example,
schematically depicted in Figure~\ref{fig:ciliary-cicle}, is the
ciliary stroke cycle. 

Our model swimmer mimics this behavior by alternating propulsive
phases (top side of the square stroke in
Figure~\ref{fig:stick-square-stroke}, i.e., $\xi_1$ varying from zero
to one while $\xi_2$ is close to one) with a recovery phases (bottom
side of the square stroke in Figure~\ref{fig:stick-square-stroke},
i.e., $\xi_1$ varying from one to zero while $\xi_2$ is close to
zero).

In order to compare the stick and donut swimmer with the three sphere
swimmer, we set the radius $R_0$ of the swimmer such that the total
volume of the stick and donut swimmer is equal to the volume of the
three sphere swimmer with radius $.05$. This gives a radius $R_0$ of
approximately $.034mm$.

\begin{figure}[!htb]
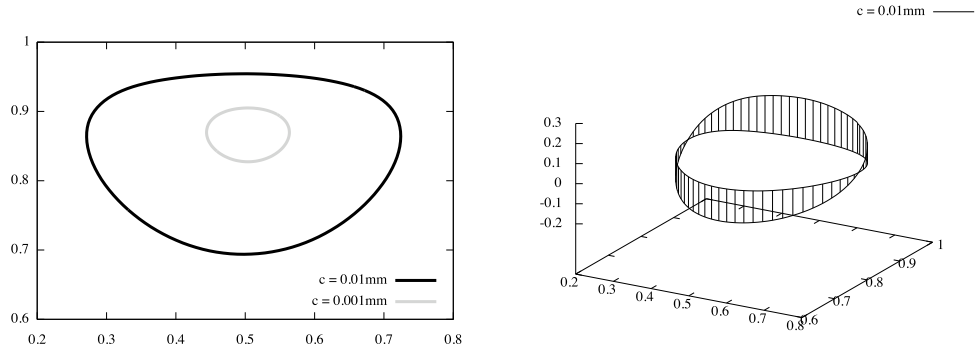

  \centering
  \includegraphics[width=.49\textwidth]{\Figpath{stick-shape}}
  \hfill
  \includegraphics[width=.49\textwidth]{\Figpath{stick-phi-dot}}
  \caption{Path in shape space describing optimal strokes for stick
    and donut swimmer with base
    radius of $.034mm$, swimming $.01mm$ and $.001mm$ in 1
    second, without imposing initial shape (left) and propulsion
    diagram for the bigger stroke (right).}
  \label{fig:stick-optimal-stroke}
\end{figure}

Figure~\ref{fig:stick-optimal-stroke} shows two optimal strokes
found with our method which lead to displacements in $1$s of
$.01mm$ and $.001mm$. In the first case, the energy consumption
is around $.126pJ$, while in the second case it is $.010pJ$.
This should be compared with the energy consumed by an
``equivalent'' three sphere swimmer, namely, a swimmer with the
same volume as this one, swimming at the same average velocity.
The three sphere swimmer energy consumption to swim  $.01mm$ in
one second is $.183pJ$, while its energy consumption to swim
$.001mm$ in one second is $.018pJ$.

In the long distance, the stick and donut swimmer is about 45\% more
efficient than the optimal three sphere swimmer, while for the shorter
distance, the increase in efficiency is about 75\%.

A collection of animations referring to the optimal strokes of the
swimmers presented in this work
(Figures~\ref{fig:3sphere-optimal-stroke}
and~\ref{fig:stick-optimal-stroke}) can be viewed on-line on the home
page of the corresponding author.\cite{Heltai--a}

\section*{Acknowledgement}
The first author benefited from the support of the
``Chair Mathematical modelling and numerical simulation,
  F-EADS -- Ecole Polytechnique -- INRIA -- F-X''.

\bibliographystyle{elsarticle-num}
\bibliography{AlougesDeSimoneHeltai}

\end{document}